\documentclass[a4paper,12pt]{amsart}

\usepackage{color}
\usepackage{epic}
\usepackage{eepic}
\usepackage{graphicx}
\usepackage{bbm}
\usepackage{url}
\usepackage{paralist}
\usepackage[tight,hang,TABTOPCAP]{subfigure}
\usepackage[small,it]{caption2}
\usepackage{tabularx}

\newtheoremstyle{plainCiteNumber}{}{}{\normalfont\itshape}{}
                {\bfseries}
                {. }{ }{\thmname{{\bfseries #1}}\thmnote{\bfseries\ #3}}

\theoremstyle{plainCiteNumber}
\newtheorem*{theorem*}{Theorem}
\newtheorem*{corollary*}{Corollary}

\theoremstyle{plain}
\newtheorem{theorem}{Theorem}[section]
\newtheorem{lemma}[theorem]{Lemma}
\newtheorem{proposition}[theorem]{Proposition}

\newtheorem*{question}{Question}
\theoremstyle{definition}
\newtheorem*{definition}{Definition}

\theoremstyle{remark}

\newcommand{\R}{\mathbbm{R}}
\newcommand{\Z}{\mathbbm{Z}}
\renewcommand{\r}{\mathbf{r}}
\renewcommand{\c}{\mathbf{c}}
\newcommand{\transport}{{\operatorname{T}_{\r\c}}}
\newcommand{\cell}[6]{$\operatorname{\mathbf
    Z}^{#1,#2,#3}_{#4,#5,#6}$}
\newcommand{\cellquer}[6]{$\operatorname{\mathbf{\overline Z}}^{#1,#2,#3}_{#4,#5,#6}$}
\newcommand{\conv}{\operatorname{conv}}

\newcommand{\subsub}[1]{\subsubsection*{#1}}

\DeclareMathOperator{\ZZ}{\overline{\mathbf Z}}

\DeclareMathOperator{\vv}{\mathbf v}

\title[Gr\"obner Bases for Transportation Polytopes]{Quadratic
  Gr\"obner Bases for Smooth $\mathbf{3 \times 3}$ Transportation
  Polytopes}

\author{Christian Haase}

\author{Andreas Paffenholz}

\address{Fachbereich Mathematik und Informatik \\
  Freie Universit\"at Berlin}

\email{\{Christian.Haase,Andreas.Paffenholz\}@Math.FU-Berlin.de}

\thanks{Both authors were supported by Emmy Noether grant HA 4383/1 of
  the German Research Foundation (DFG)}

\begin{document}

\begin{abstract}
  The toric ideals of $3 \times 3$ transportation polytopes
  $\transport$ are quadratically generated. The only exception is the
  Birkhoff polytope $B_3$.

  If $\transport$ is not a multiple of $B_3$, these ideals even have
  squarefree quadratic initial ideals. This class contains all smooth
  $3 \times 3$ transportation polytopes.
\end{abstract}

\maketitle

\section{Introduction}
\subsection{Motivation$^1$}
\stepcounter{footnote}
\footnotetext{This motivation is quoted verbatim from~\cite{CottonwoodRoom}.}
A lattice polytope $P \subset \R^d$ defines an ample line bundle $L_P$
on a projective toric variety $X_P$. (See, e.g.,
\cite[\S 3.4]{FultonToric}.) If $X_P$ is smooth (the normal fan of $P$ is
unimodular), then $L_P$ is very ample, and provides an embedding $X_P
\hookrightarrow \mathbb{P}^{r-1}$, where $r = \# (P \cap \Z^d)$. So we
can think of $X_P$ as canonically sitting in projective space.
The following question \cite[Conjecture 2.9]{SturmfelsToricEquations}
about the defining equations of $X_P \subset \mathbb{P}^{r-1}$ has
been around for quite a while, but its origins are hard to track
(cf.~\cite{CottonwoodRoom}).
\begin{question} \label{q:quadratic}
  Let $P$ be a lattice polytope whose corresponding projective toric
  variety is smooth. Is the defining ideal $I_P$ generated by
  quadratics?
\end{question}
There are   two variations of this question (which  are  of strictly
increasing strength).
\begin{itemize}
  \item Is the homogeneous coordinate ring
  $\mathbbm{k}[x_1, \ldots, x_r]/I_P$ Koszul?
  \item Does $I_P$ have a quadratic Gr\"obner basis?
\end{itemize}
The last version has a combinatorial interpretation. It asks for the
existence of very special, ``quadratic'' triangulations of $P$, see
\S\ref{sec:background} below.

\subsection{Results}
Simple transportation polytopes provide a large family of smooth
polytopes. Yet, the $3 \times 3$ Birkhoff polytope $B_3$ is a
non-simple transportation polytope whose ideal is not generated
by quadratic polynomials. In this note, we show that in the $3 \times
3$ case, this is the only example.
In Section~\ref{sec:quadratic-generation}, we show that $B_3$ is the
only ($3 \times 3)$-transportation polytope whose ideal is not
quadratically generated.
\begin{proposition} \label{prop:generation}
  If $\transport \neq B_3$, then $I_{\transport}$ is quadratically
  generated.
\end{proposition}

If $P$ is a $3 \times 3$ transportation polytope which is not a
multiple of $B_3$, we can show in
Section~\ref{sec:quadr-groebner-bases} that these ideals even have
quadratic Gr\"obner bases. This class contains all smooth $3
\times 3$ transportation polytopes. 
\newcommand{\textofcorquadratic}{If  $\transport$ is not a multiple of
  $B_3$, then $I_{\transport}$  has a  squarefree quadratic initial
  ideal.}
\begin{theorem} \label{cor:GB}
\textofcorquadratic
\end{theorem}

Using different methods, Lindsay Piechnik and the first author showed
that (among other polytopes) even multiples of $B_3$ have quadratic
triangulations. We believe that odd multiples $\ge 3$ allow quadratic
triangulations as well.

\subsection{Background}\label{sec:background}
\subsub{Transportation Polytopes}
Let two vectors $\c=(c_1,\ldots, c_n)\in\Z_{>0}^n$ and
$\r=(r_1,\ldots, r_m)\in\Z_{>0}^m$ with
$\sum_{i=1}^nc_i=\sum_{i=1}^mr_i=:s$ be given. The corresponding
{$(m\times n)$-transportation polytope} $\transport$ is the set
of all non-negative $(m \times n)$-matrices $A=(a_{ij})_{ij}$ satisfying
\begin{align*}
  \sum_{i=1}^ma_{ik}=c_k\qquad\text{and}\qquad\sum_{j=1}^na_{lj}=r_l
\end{align*}
for $1 \le k \le n$,  $1 \le l \le m$.  This is a
bounded convex polytope with integral vertices (a lattice polytope for
short) in $\R^{mn}$.  We  number the  coordinates of $\R^{mn}$ by
$a_{ij}$ for  $1\le i\le m$ and  $1\le j\le n$.
The upper $((m-1)\times (n-1))$-minor of a matrix $A$ in the polytope
determines all other  entries. Hence, the dimension of  $\transport$
is at most $(m-1)(n-1)$. On the other hand, $a_{ij}=r_ic_j/s$
determines an interior point, so that the dimension is exactly
$(m-1)(n-1)$. In what follows, we focus on the case $m=n=3$.
%%%%%%%
\subsub{Toric Ideals}
Let $P \subset \R^d$ be a lattice polytope. The point configuration
$\mathcal A = P \cap \Z^d = \{ \mathbf a_1,
\ldots, \mathbf a_r \}$ defines a ring homomorphism
$$
\begin{array}{ccc}
  \mathbbm k[x_1, \ldots, x_r] & \longrightarrow & \mathbbm
  k[t_0^{\pm 1}, \ldots t_d^{\pm 1}] \\[1mm]
  x_i & \longmapsto & t_0 \, \mathbf t^{\mathbf a_i} := t_0 \,
  t_1^{a_{1i}} \cdot \ldots \cdot t_d^{a_{d i}}.
\end{array}
$$
Its kernel is the homogenous  ideal $$I_P = \langle \mathbf{x^u-x^v} \
: \ \sum u_i \mathbf a_i = \sum v_i \mathbf a_i \ , \  \sum u_i = \sum
v_i \rangle .$$ This ideal is called the toric ideal associated to $P$
(see~\cite[\S4]{SturmfelsGBCP}).
%%%%%%%
\subsub{The Birkhoff Polytope}
The simplest
$(3 \times 3)$-transportation polytope is the Birkhoff polytope $B_3$
of doubly stochastic matrices, given by $\r = \c = (1,1,1)$. 
% $B_3$ is not simple, and hence, not smooth.
The lattice points in $B_3$ are the six permutation matrices
$A_\sigma$ for $\sigma \in S_3$. If we denote the corresponding
variables by $x_\sigma$, the toric ideal $I_{B_3}$ is the principal
ideal $\langle x_{123}x_{231}x_{312} - x_{132}x_{213}x_{321}
\rangle$. So $I_{B_3}$ is not quadratically generated. $I_{B_3}$ has
two initial ideals, $\langle x_{123}x_{231}x_{312} \rangle$, and
$\langle x_{132}x_{213}x_{321} \rangle$.
Geometrically, this correponds to the fact that $B_3 \cap \Z^9$ is a
circuit, i.e., a minimal affinely dependent set. $B_3$ is the convex
hull of the triangle of even permutation matrices together with the
triangle of odd permutation matrices. The two triangles meet in their
barycenters.
\begin{equation}
  \label{eq:B3} \tag{$\ast$}
  \left[
    \begin{smallmatrix}
      1&&\\&1&\\&&1
    \end{smallmatrix}
  \right]
  +
  \left[
    \begin{smallmatrix}
      &1&\\&&1\\1&&
    \end{smallmatrix}
  \right]
  +
  \left[
    \begin{smallmatrix}
      &&1\\1&&\\&1&
    \end{smallmatrix}
  \right]
  =
  \left[
    \begin{smallmatrix}
      1&&\\&&1\\&1&
    \end{smallmatrix}
  \right]
  +
  \left[
    \begin{smallmatrix}
      &1&\\1&&\\&&1
    \end{smallmatrix}
  \right]
  +
  \left[
    \begin{smallmatrix}
      &&1\\&1&\\1&&
    \end{smallmatrix}
  \right]
\end{equation}
This (up to scaling) unique affine relation yields the equation
generating $I_{B_3}$. 
%%%%%%%
\subsub{Smooth Polytopes} For a lattice polytope $P$, the set of zeros
in $\mathbbm P^{r-1}$ of $I_P$ is the toric variety $X_P$. This
variety is smooth if and only if the edge directions at every vertex
of $P$ form a lattice basis. Equivalently, $X_P$ is smooth if and only
if the normal fan of $P$ is unimodular (See~\cite[\S
2.1]{FultonToric}). In this case we call $P$ a smooth polytope. In
particular, smooth polytopes are simple: every vertex belongs to
dimension many facets. (So, the Birkhoff polytope is not smooth.)
\begin{lemma}\label{lemma:smoothTransport}
  For a transportation polytope $\transport$, the following are
  equivalent.
  \begin{compactenum} 
  \item $X_{\transport}$ is smooth.
  \item $\transport$ is smooth.
  \item $\transport$ is simple.
  \item $\sum_{i \in I}r_i \neq \sum_{j \in J}c_j$ for all non-trivial
    sets of indices $I \subset [m]$, $J \subset [n]$.
  \end{compactenum}
\end{lemma}
We have not found a proof in the literature. For completeness, we
include one here. (Compare the discussion for general flow polytopes
in~\cite{BSdLV}. Lemma~\ref{lemma:smoothTransport} says that in our
case, topes and chambers agree.) 
\begin{proof}
  \underline{(1) $\Leftrightarrow$ (2)} by~\cite[\S~2.1]{FultonToric}.
  The implication \underline{(2) $\Rightarrow$ (3)} is valid for all
  lattice polytopes. The converse, \underline{(3) $\Rightarrow$ (2)}
  follows from the fact that transportation polytopes arise from a
  totally unimodular matrix~\cite[\S 19]{Schrijver}.

  \underline{(4) $\Rightarrow$ (3)}:
  Suppose that $\transport$ has a vertex $A$ that belongs to $\ge
  (m-1)(n-1)+1$ facets. Then $A$ has at least that many zero
  entries. Thus, the bipartite graph given by the non-zero entries has
  $n+m$ vertices and $\le n+m-2$ edges. So this graph cannot be
  connected. Take for $I$ and $J$ the color classes of one component
  of this graph.

  For \underline{(3) $\Rightarrow$ (4)} we need some preliminary
  observations. We use the criterion that an inequality $a_{ij} \ge 0$
  defines a facet of $\transport$ if and only if there is an $A \in
  \transport$ such that $a_{ij} = 0$ and with all other entries
  positive. 

  Now, suppose we are given $I \subset [m]$ and $J \subset [n]$ with
  $\sum_{i \in I}r_i = \sum_{j \in J}c_j$. 
  Build a matrix $A \in \transport$ from a vertex $A'$ of the $I
  \times J$ transportation polytope, and a vertex $A''$ of the $I^c
  \times J^c$ transportation polytope. We abbreviate $m'=|I|$,
  $m''=|I^c|$, $n'=|J|$, and $n''=|J^c|$.
  \begin{center}
    \begin{tabular}{c|c|c|}
      \multicolumn{1}{c}{} & \multicolumn{1}{c}{$I$} &
      \multicolumn{1}{c}{$I^c$} \\[1mm]
      \cline{2-3}
      $J$ \raisebox{6mm}{\mbox{}} & $A'$ & $0$ \\[1mm]
      \cline{2-3}
      $J^c$ \raisebox{6mm}{\mbox{}} & $0$ & $A''$ \\[1mm]
      \cline{2-3}
    \end{tabular}
  \end{center}
  \begin{lemma}
    The inequalities $a_{ij} \ge 0$ for $(i,j) \in I \times J^c \cup I^c
    \times J$ define facets of $\transport$.
  \end{lemma}
  \begin{proof}
    Say, $(i,j) \in I \times J^c$. Start from all positive $A'$ and
    $A''$. Add $\varepsilon m''n'$ to all $I \times J^c$ entries $\neq
    (i,j)$, and $\varepsilon (m'n''-1)$ to all $I^c \times J$
    entries. Now modify $A'$ and $A''$ in order to obtain the old row
    and column sums. This amounts to finding points in two
    (non-integral) transportation polytopes. For small enough
    $\varepsilon$, the resulting matrix will have positive entries
    away from $(i,j)$.
  \end{proof}
  \begin{lemma}
    If the inequality $a_{ij} \ge 0$ ($(i,j) \in I \times J$) defines a
    facet of the $I \times J$ transportation polytope, then it also
    defines a facet of $\transport$.
  \end{lemma}
  \begin{proof}
    Let $A'$ be a matrix whose only zero entry is $(i,j)$, and let
    $A''$ be all positive. As before, we can subtract suitable
    constants from $A'$ and $A''$, and find all positive matrices to
    couterbalance row and column sums.
  \end{proof}
  To wrap it up, if $A'$ and $A''$ are vertices of their transportation
  polytopes, the block matrix $A$ belongs to at least $(m-1)(n-1)+1$
  facets. Hence, $\transport$ is not simple.
\end{proof}
%%%%%%%
\subsub{Triangulations}
In order to show that a toric ideal has a quadratic Gr\"obner basis,
we use the connection to regular triangulations as outlined
in~\cite[\S8]{SturmfelsGBCP}. A subset $F \subseteq P \cap \Z^d$ is a
face of a triangulation of $P$ if $\conv(F)$ is a simplex of the
triangulation; otherwise $F$ is said to be a non-face. Observe that
every superset of a non-face is a non-face.
\begin{definition}
  A regular unimodular triangulation whose minimal non-faces have two
  elements is called a quadratic triangulation.\footnote{Simplicial
    complexes with this non-face property appear in the literature
    under the names of flag- or clique-complexes.}
\end{definition}
The following characterization is a conglomerate of Corollaries~8.4
and 8.9 in~\cite{SturmfelsGBCP}.
\begin{theorem} \label{thm:GBtriang}
  The defining ideal $I_P$ of the projective toric variety $X_P
  \subset \mathbb{P}^{r-1}$ has a squarefree initial ideal
  if and only if $P$ has a regular unimodular triangulation.

  In that case, the corresponding initial ideal is the Stanley-Reisner
  ideal of the triangulation: $\operatorname{in}(I_P) = \langle x^F \
  \vert \ F \text{ minimal non-face } \rangle$ .
\end{theorem}
Here, we abbreviate $x^F := \prod_{i \in F} x_i$.
In the example of the Birkhoff polytope, there are two (isomorphic)
triangulations of $B_3$. They are both regular and unimodular. In one
of them, the triangle of even permutation matrices is the minimal
non-face, in the other one, the triangle of odd permutation matrices
is the minimal non-face.

Using this
correspondence, Theorem~\ref{cor:GB} follows from the following
theorem which is what we really prove in
Section~\ref{sec:quadr-groebner-bases}.
%%%
\newcommand{\textofthmquadratic}{If $\transport$ is  not a multiple of
  $B_3$, then $\transport$ has a quadratic triangulation.}
%%%
\begin{theorem} \label{thm:GB}
  \textofthmquadratic
\end{theorem}
%%%%%%%
\subsub{Paco's Lemma}  A tool  we   use  in  both proofs  are  pulling
refinements of hyperplane  subdivisions.  Let  $P \subset \R^d$   be a
lattice polytope. As before, order the lattice points $P \cap \Z^d =
\mathbf a_1, \ldots, \mathbf a_r$, and the corresponding  variables
$x_1 \prec \ldots \prec x_r$. Then, the  reverse lexicographic term
order  yields  a {pulling triangulation\/} of $P$. These  pulling
triangulations have a nice recursive structure: 
the maximal faces are joins of $\mathbf a_1$ with faces of the
pulling triangulations of those facets of $P$ that do not contain
$\mathbf a_1$.

We say that a lattice polytope $P$ has facet  width $1$ if for each of
its facets, $P$ lies between the hyperplane  spanned by this facet and
the next parallel lattice hyperplane.
\begin{proposition}[Paco's Lemma~\cite{pacoPersonal,
    OhsugiHibiCompressed, SullivantCompressed}] \label{prop:paco}
  The lattice polytope $P$ has facet width $1$ if and only if every
  pulling triangulation of $P$ is (regular and) unimodular.
\end{proposition}

%%%%%%%%%%%%%%%%%%%%%%%%%%%%%%%%%%%%%%%%%%%%%%%
\section{Quadratic Generation} \label{sec:quadratic-generation}

The main tools in the proof of Proposition~\ref{prop:generation} are a
hyperplane subdivision and matrix addition. We will first exhibit a
Gr\"obner basis which consists of quadratic and cubic binomials. Then
we go on to show that the cubic elements can be expressed using
quadratic members of the ideal. The resulting quadratic generating set
will usually fail to be a Gr\"obner basis.

A transportation polytope
has a canonical regular subdivision into polytopes of facet width $1$.
We  slice  $\transport$  along    the  hyperplanes  $a_{ij}=k$.     By
Proposition~\ref{prop:paco},     every    pulling  refinement of  this
subdivision will be a regular unimodular triangulation. A non-face $F$
of such  a  triangulation either contains a   pair  of matrices  which
differ by $\ge  2$  in one entry  (a  minimal non-face of  cardinality
$2$), or   all  of $F$ belongs  to   the same cell of  the  hyperplane
subdivision.

The ideal  $I_\transport$ is generated by  a  Gr\"obner basis which is
parameterized by the   minimal non-faces of  the given  triangulation.
(And  the  degree of    a generator equals    the cardinality  of  the
corresponding non-face.)  So  we need  to   analyze the  cells of  the
hyperplane subdivision.  They have the form
$$
  \operatorname{\mathbf Z}_{\r\c}(K) = \left\{ A \in \transport \ \vert \
    k_{ij} \le a_{ij} \le k_{ij} + 1 \right\}
  %\text{ for } i=1,\ldots,m \ j=1,\ldots,n \right\}
$$
for some matrix $K$ with row  sums $\r'$ and  column sums $\c'$. After
translation we get
$$
  \operatorname{\mathbf Z}_{\r\c}(K)-K = \operatorname{\mathbf
    Z}_{\r-\r',\c-\c'}(0) =: \operatorname{\mathbf Z}^{\r-\r'}_{\c-\c'}.
$$
%In the ($3 \times 3)$-case, there are only four types of such cells.
In order to obtain a full-dimensional cell, $\r-\r'$ and $\c-\c'$ must
have coefficients $1$ or  $2$. So, up to  symmetry, in the ($3  \times
3)$-case     there are only   four   types    of such  cells,   namely
\cell{1}{1}{1}{1}{1}{1},                      \cell{1}{1}{2}{1}{1}{2},
\cell{1}{2}{2}{1}{2}{2}, and \cell{2}{2}{2}{2}{2}{2}.       In   fact,
\cell{1}{1}{2}{1}{1}{2}   and  \cell{1}{2}{2}{1}{2}{2} are  unimodular
simplices,      and           \cell{1}{1}{1}{1}{1}{1}$=B_3$        and
\cell{2}{2}{2}{2}{2}{2} are isomorphic as lattice polytopes.%
%\footnote{The isomorphism is given by $A \mapsto
%  \mathbbm{1}-A$.}

To summarize,    $I_\transport$ is  generated by   quadratic binomials
together with cubic binomials that  correspond to affine relations \`a
la~\eqref{eq:B3}.

Now     let   us  assume    $\transport       \neq B_3$,  and,    say,
$\operatorname{Z}_{\r\c}(K)  \cong$  \cell{1}{1}{1}{1}{1}{1} is a cell
in   $\transport$ giving  rise  to such   a  cubic equation.   Because
$\transport \neq  B_3$,  at  least  one of  the nine    adjacent cells
$\operatorname{Z}_{\r\c}(K-E_{ij})$ has  to be in  $\transport$, where
$E_{ij}$ is the $(ij)^\mathrm{th}$  unit vector. After translation, we
are given the relation~\eqref{eq:B3}, and we know that (for $i=j=1$)
$$
  \left[
    \begin{smallmatrix}
      -1&1&1 \\ 1 \\ 1
    \end{smallmatrix}
  \right]
  \quad \in \quad \transport-K .
$$
But then, we can use the two quadratic relations
$$
\begin{array}{ccccccc}
  \left[
    \begin{smallmatrix}
      -1&1&1 \\ 1 \\ 1
    \end{smallmatrix}
  \right]
  &+&
  \left[
    \begin{smallmatrix}
      1&&\\&1&\\&&1
    \end{smallmatrix}
  \right]
  &=&
  \left[
    \begin{smallmatrix}
      &1&\\1&&\\&&1
    \end{smallmatrix}
  \right]
  &+&
  \left[
    \begin{smallmatrix}
      &&1\\&1&\\1&&
    \end{smallmatrix}
  \right]
  \\[5mm]
  \left[
    \begin{smallmatrix}
      -1&1&1 \\ 1 \\ 1
    \end{smallmatrix}
  \right]
  &+&
  \left[
    \begin{smallmatrix}
      1&&\\&&1\\&1&
    \end{smallmatrix}
  \right]
  &=&
  \left[
    \begin{smallmatrix}
      &&1\\1&&\\&1&
    \end{smallmatrix}
  \right]
  &+&
  \left[
    \begin{smallmatrix}
      &1&\\&&1\\1&&
    \end{smallmatrix}
  \right]
\end{array}
$$
to  generate~\eqref{eq:B3}.   This      completes the     proof     of
Proposition~\ref{prop:generation}.  \qed
%%%%%%%%%%%%%%%%%%%%%%%%%%%%%%%%%%%
\section{Quadratic Gr\"obner Bases} \label{sec:quadr-groebner-bases}
In the previous section we   have seen that  all  toric   ideals
associated to  transportation polytopes  $\transport \neq B_3$ are
quadratically generated. Now we head for the stronger result stated in
Theorem~\ref{cor:GB}.

We again start by subdividing  the polytope into cells by intersection
with hyperplanes  of the form  $a_{ij}=k$,  but this  time we choose a
coarser subdivision to   avoid   cells  isomorphic to    the  Birkhoff
polytope. We show that  we can achieve  this by taking all hyperplanes
$a_{ij}=k$ except for $(i,j)=(1,1)$ and $(i,j)=(2,1)$.  In a second
step we do a pulling triangulation along a vertex order given by a
(globally fixed) linear functional.
The analysis of the cells was done using the software package {\tt
  polymake}~\cite{polymake}.
%%%%%%%%%%%
\subsection{Hyperplane Subdivision}
Let $\transport$ be a transportation polytope with row and column sums
$\r$  and $\c$, which is not  a multiple of  $B_3$.  We order the rows
and columns in such a way that $r_1\ge  r_2\ge r_3$ and $c_1\ge c_2\ge
c_3$. We can assume that $c_1  \ge r_1$, and  thus $c_1 > r_3$ because
$\transport$ is not a multiple of $B_3$.

As  before, we subdivide $\transport$ into  cells  by cutting with the
hyperplanes  $\{a_{ij}=k_{ij}\}$ except that we  omit the $(1,1)$- and
the $(2,1)$-entries. Hence, our cells are of the form
\begin{align*}
  \ZZ_{\r\c}(K) = \left\{ A\in\transport
    \;\rule[-.8\baselineskip]{.5pt}{2\baselineskip}\;
    {\genfrac{}{}{0pt}{}
      {k_{ij} \le a_{ij} \le k_{ij} + 1\hfill }
      {\text{for }(i,j)=(3,1)\text{ and }1\le i\le 3,\; 2\le j\le 3}} 
  \right\}
\end{align*}
where     $K=(k_{ij})_{ij}$ is     a    $(3\times   3)$-matrix    with
$k_{11},k_{21}:=0$.
Similar to the previous section, we subtract $K$ from $\transport$,
arriving at cells of the form
\begin{align*}
  \ZZ_{\r\c}(K)-K = \ZZ_{\r-\r',\c-\c'}(0) .
\end{align*}
Again, we get $r_3-r_3', c_2-c_2', c_3-c_3' \in \{1,2\}$.
Also, we have
\begin{align*}
  c_1-c_1' = c_1 - k_{31} > r_3 - k_{31} \ge r_3 - r_3'.
\end{align*}

The projection to the $a_{11}$-$a_{21}$-plane is described by the
inequalities $r_i-r_i'-2 \le a_{i1} \le r_i-r_i'$ for $i=1,2$ and
$c_1-c_1'-1 \le a_{11}+a_{21} \le c_1-c_1'$. Thus, if $r_1,r_2 \ge
2$, there are four cases (cf. Figure~\ref{fig:projection} on the
left). If $r_1-r_1' = 1$ or $r_2-r_2' = 1$, there are three and three
more cases (cf. Figure~\ref{fig:projection} on the right). We cannot
have $r_1-r_1' = r_2-r_2' = 1$ because $c_1-c_1'>r_3-r_3'$.

\newcommand{\qq}{\scriptstyle}
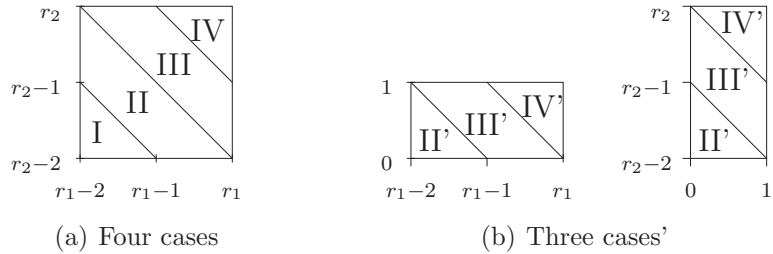
\begin{figure}[htbp]
  \centering
\begingroup\makeatletter\ifx\SetFigFont\undefined%
\gdef\SetFigFont#1#2#3#4#5{%
  \reset@font\fontsize{#1}{#2pt}%
  \fontfamily{#3}\fontseries{#4}\fontshape{#5}%
  \selectfont}%
\fi\endgroup%
\subfigure[Four cases]{
\setlength{\unitlength}{0.00043745in}
{\renewcommand{\dashlinestretch}{30}
\begin{picture}(2055,2589)(0,-10)
\path(240,2382)(240,582)(2040,582)
	(2040,2382)(240,2382)
\path(240,1482)(1140,582)
\path(240,2382)(2040,582)
\path(1140,2382)(2040,1482)
\path(240,627)(240,537)
\path(1140,627)(1140,537)
\path(2040,627)(2040,537)
\path(195,2382)(285,2382)
\path(195,1482)(285,1482)
\path(195,582)(285,582)
\put(1140,132){\makebox(0,0)[b]{\smash{{\SetFigFont{12}{14.4}{\rmdefault}{\mddefault}{\updefault}$\qq r_1-1$}}}}
\put(15,1347){\makebox(0,0)[rb]{\smash{{\SetFigFont{12}{14.4}{\rmdefault}{\mddefault}{\updefault}$\qq r_2-1$}}}}
\put(15,2247){\makebox(0,0)[rb]{\smash{{\SetFigFont{12}{14.4}{\rmdefault}{\mddefault}{\updefault}$\qq r_2$}}}}
\put(15,447){\makebox(0,0)[rb]{\smash{{\SetFigFont{12}{14.4}{\rmdefault}{\mddefault}{\updefault}$\qq r_2-2$}}}}
\put(375,762){\makebox(0,0)[lb]{\smash{{\SetFigFont{12}{14.4}{\rmdefault}{\mddefault}{\updefault}I}}}}
\put(780,1122){\makebox(0,0)[lb]{\smash{{\SetFigFont{12}{14.4}{\rmdefault}{\mddefault}{\updefault}II}}}}
\put(1140,1572){\makebox(0,0)[lb]{\smash{{\SetFigFont{12}{14.4}{\rmdefault}{\mddefault}{\updefault}III}}}}
\put(1545,1977){\makebox(0,0)[lb]{\smash{{\SetFigFont{12}{14.4}{\rmdefault}{\mddefault}{\updefault}IV}}}}
\put(240,132){\makebox(0,0)[b]{\smash{{\SetFigFont{12}{14.4}{\rmdefault}{\mddefault}{\updefault}$\qq r_1-2$}}}}
\put(2040,132){\makebox(0,0)[b]{\smash{{\SetFigFont{12}{14.4}{\rmdefault}{\mddefault}{\updefault}$\qq r_1$}}}}
\end{picture}
}}
%  \qquad
\hspace{.1\textwidth}
\begingroup\makeatletter\ifx\SetFigFont\undefined%
\gdef\SetFigFont#1#2#3#4#5{%
  \reset@font\fontsize{#1}{#2pt}%
  \fontfamily{#3}\fontseries{#4}\fontshape{#5}%
  \selectfont}%
\fi\endgroup%
\subfigure[Three cases']{
\setlength{\unitlength}{0.00043745in}
{\renewcommand{\dashlinestretch}{30}
\begin{picture}(2083,1689)(0,-10)
\path(240,1482)(1140,582)
\path(1140,1482)(2040,582)
\path(195,1482)(285,1482)
\path(195,582)(285,582)
\path(240,1482)(240,582)(2040,582)
	(2040,1482)(240,1482)
\path(240,582)(240,492)
\path(1140,582)(1140,492)
\path(2040,627)(2040,537)
\put(240,132){\makebox(0,0)[b]{\smash{{\SetFigFont{12}{14.4}{\rmdefault}{\mddefault}{\updefault}$\qq r_1-2$}}}}
\put(1140,132){\makebox(0,0)[b]{\smash{{\SetFigFont{12}{14.4}{\rmdefault}{\mddefault}{\updefault}$\qq r_1-1$}}}}
\put(2040,132){\makebox(0,0)[b]{\smash{{\SetFigFont{12}{14.4}{\rmdefault}{\mddefault}{\updefault}$\qq r_1$}}}}
\put(15,447){\makebox(0,0)[rb]{\smash{{\SetFigFont{12}{14.4}{\rmdefault}{\mddefault}{\updefault}$\qq 0$}}}}
\put(15,1347){\makebox(0,0)[rb]{\smash{{\SetFigFont{12}{14.4}{\rmdefault}{\mddefault}{\updefault}$\qq 1$}}}}
\put(330,672){\makebox(0,0)[lb]{\smash{{\SetFigFont{12}{14.4}{\rmdefault}{\mddefault}{\updefault}II'}}}}
\put(870,852){\makebox(0,0)[lb]{\smash{{\SetFigFont{12}{14.4}{\rmdefault}{\mddefault}{\updefault}III'}}}}
\put(1545,1077){\makebox(0,0)[lb]{\smash{{\SetFigFont{12}{14.4}{\rmdefault}{\mddefault}{\updefault}IV'}}}}
\end{picture}
}
\quad \quad
\setlength{\unitlength}{0.00043745in}
\begingroup\makeatletter\ifx\SetFigFont\undefined%
\gdef\SetFigFont#1#2#3#4#5{%
  \reset@font\fontsize{#1}{#2pt}%
  \fontfamily{#3}\fontseries{#4}\fontshape{#5}%
  \selectfont}%
\fi\endgroup%
{\renewcommand{\dashlinestretch}{30}
\begin{picture}(1155,2589)(0,-10)
\path(240,1482)(1140,582)
\path(240,2382)(240,582)(1140,582)
	(1140,2382)(240,2382)
\path(240,582)(240,492)
\path(1140,582)(1140,492)
\path(240,2382)(1140,1482)
\path(150,1482)(240,1482)
\path(150,582)(240,582)
\path(150,2382)(240,2382)
\put(240,132){\makebox(0,0)[b]{\smash{{\SetFigFont{12}{14.4}{\rmdefault}{\mddefault}{\updefault}$\qq 0$}}}}
\put(15,447){\makebox(0,0)[rb]{\smash{{\SetFigFont{12}{14.4}{\rmdefault}{\mddefault}{\updefault}$\qq r_2-2$}}}}
\put(15,1347){\makebox(0,0)[rb]{\smash{{\SetFigFont{12}{14.4}{\rmdefault}{\mddefault}{\updefault}$\qq r_2-1$}}}}
\put(330,672){\makebox(0,0)[lb]{\smash{{\SetFigFont{12}{14.4}{\rmdefault}{\mddefault}{\updefault}II'}}}}
\put(420,1392){\makebox(0,0)[lb]{\smash{{\SetFigFont{12}{14.4}{\rmdefault}{\mddefault}{\updefault}III'}}}}
\put(600,2067){\makebox(0,0)[lb]{\smash{{\SetFigFont{12}{14.4}{\rmdefault}{\mddefault}{\updefault}IV'}}}}
\put(15,2247){\makebox(0,0)[rb]{\smash{{\SetFigFont{12}{14.4}{\rmdefault}{\mddefault}{\updefault}$\qq r_2$}}}}
\put(1140,132){\makebox(0,0)[b]{\smash{{\SetFigFont{12}{14.4}{\rmdefault}{\mddefault}{\updefault}$\qq 1$}}}}
\end{picture}
}

}
%  \caption{Four Cases. Three Cases'.}
\caption{The projection onto the $a_{11}$-$a_{21}$-plane}
  \label{fig:projection}
\end{figure}

If we  subtract the lower bounds for  $a_{11}$ and $a_{21}$, we obtain
the following $20$  translation classes of  cells  in the subdivision.
(We list them in the form $\mathbf{(r-r')\ (c-c')}$.)
\begin{itemize}
\item[I] $(2,2,1)(1,2,2)$
\item[II] $(2,2,1)(2,2,1)$ , $(2,2,1)(2,1,2)$ , $(2,2,2)(2,2,2)$
\item[III] $(2,2,1)(3,1,1)$ , $(2,2,2)(3,2,1)$ , $(2,2,2)(3,1,2)$
\item[IV] $(1,1,2)(2,1,1)$
\item[II'] $(2,1,1)(1,2,1)$ , $(2,1,1)(1,1,2)$ , $(2,1,2)(1,2,2)$ ,
%  and symmetrically, \\
  \\
  $(1,2,1)(1,2,1)$ , $(1,2,1)(1,1,2)$ , $(1,2,2)(1,2,2)$
\item[III'] $(2,1,1)(2,1,1)$ , $(2,1,2)(2,2,1)$ , $(2,1,2)(2,1,2)$ , 
%  and symmetrically, \\
  \\
  $(1,2,1)(2,1,1)$ , $(1,2,2)(2,2,1)$ , $(1,2,2)(2,1,2)$
\item[IV'] same as IV.
\end{itemize}

\subsection{Triangulating the Cells}
According to 
%Paco's 
Lemma~\ref{prop:paco}, any pulling refinement of
our cell decomposition will be unimodular. The subtle part is to
devise a pulling order so that the resulting triangulation is flag,
i.e., so that the minimal non-faces have cardinality two.

Just as before, a non-face of such  a  triangulation either contains a
non-face of  cardinality $2$, or it belongs to the same cell of the
hyperplane subdivision.
Hence, it suffices to guarantee that the induced traingulations of the
cells are flag. To achieve this, we order the vertices by decreasing
values of the linear functional
\begin{align*}
  \vv:= %[ 4, 6, -1, 3] .
  \left[
  \begin{smallmatrix}
    4&6&0\\-1&3&0\\0&0&0
  \end{smallmatrix}
  \right]
\end{align*}
The induced triangulations of the cells are the pulling triangulations
in the induced vertex ordering.

\subsubsection{Description of the Cells}
We give geometric descriptions for all possible cells. For most of
them, any triangulation is unimodular and flag. There are two
combinatorial types where we have to be careful. We list the vertices
of those cells explicitely in Table~\ref{tab:forthset} in the form
$[a_{11} \ a_{12} \ a_{21} \ a_{22}]$.

\subsub{\cellquer{1}{1}{2}{2}{1}{1}, \cellquer{1}{2}{1}{1}{1}{2},
  \cellquer{1}{2}{1}{1}{2}{1}, \cellquer{1}{2}{2}{1}{2}{2},
  \cellquer{2}{1}{1}{1}{1}{2}, \cellquer{2}{1}{1}{1}{2}{1}, 
  \cellquer{2}{1}{2}{1}{2}{2}, and \cellquer{2}{2}{1}{1}{2}{2}
}
These cells already are unimodular simplices.

\subsub{\cellquer{2}{2}{1}{2}{1}{2}, \cellquer{2}{2}{1}{2}{2}{1},
  \cellquer{2}{2}{2}{3}{1}{2}, and \cellquer{2}{2}{2}{3}{2}{1}}

These cells are pyramids over a triangular prism $\Delta_2 \times
\Delta_1$.
All six trianglulations of such a polytope are unimodular and flag.

\subsub{\cellquer{1}{2}{2}{2}{1}{2}, \cellquer{1}{2}{2}{2}{2}{1},
  \cellquer{2}{1}{2}{2}{2}{1}, and \cellquer{2}{1}{2}{2}{1}{2}}

These cells are a join of an edge and a unit square.
Both trianglulations of such a polytope are unimodular and flag.

\subsub{\cellquer{1}{2}{1}{2}{1}{1} and \cellquer{2}{1}{1}{2}{1}{1}}

These cells are isomorphic to a Birkhoff polytope $B_3$ where we have
relaxed one facet. Their vertices are listed in
Table~\ref{tab:z121-211} and ~\ref{tab:z211-211}.

\subsub{\cellquer{2}{2}{1}{3}{1}{1} and \cellquer{2}{2}{2}{2}{2}{2}}

These cells are combinatorial duals of the cells in the previous
paragraph. Their vertices are in Tables~\ref{tab:z222-222} and
\ref{tab:z221-311}.

\begin{table}[bt]
    \caption{Birkoff with one relaxed facet and its dual. We
      record the value of our linear functional in front of each vertex.}
    \label{tab:forthset}  
  \begin{minipage}[b]{.48\textwidth}
    \footnotesize
    \subtable[\cellquer{1}{2}{1}{2}{1}{1}]{
      \begin{minipage}[b]{.47\textwidth}
        \begin{alignat*}{6}
      7:&&\;\;    [1&&\;\; 0&&\;\; 0&&\;\; 1]\\[-.05cm]
      6:&&    [1&&0&&1&&1]\\[-.05cm]
      5:&&    [0&&1&&1&&0]\\[-.05cm]
      4:&&    [0&&1&&2&&0]\\[-.05cm]
      3:&&    [1&&0&&1&&0]\\[-.05cm]
      2:&&    [0&&0&&1&&1]\\[-.05cm]
     -2:&&    [0&&0&&2&&0]
        \end{alignat*}
        \label{tab:z121-211}
      \end{minipage}}
    \subtable[\cellquer{2}{1}{1}{2}{1}{1}]{
      \begin{minipage}[b]{.47\textwidth}
        \begin{alignat*}{6}
      11:&&\;\;    [2&&0&&0&&1]\\[-.05cm]
      10:&&    [1&&\;\; 1&&\;\; 0&&\;\; 0]\\[-.05cm]
      9:&&    [1&&1&&1&&0]\\[-.05cm]
      8:&&    [2&&0&&0&&0]\\[-.05cm]
      7:&&    [1&&0&&0&&1]\\[-.05cm]
      5:&&    [0&&1&&1&&0]\\[-.05cm]
      3:&&    [1&&0&&1&&0]
        \end{alignat*}
        \label{tab:z211-211}
      \end{minipage}}
  \end{minipage}
  \begin{minipage}[b]{.48\textwidth}
    \footnotesize
    \subtable[\cellquer{2}{2}{2}{2}{2}{2}]{
      \begin{minipage}[b]{.47\textwidth}
        \begin{alignat*}{6}
       13:&&\;\;   [&&1&& 1&& 0&& 1]\\[-.05cm]
       11:&&   [&&2&& 0&& 0&& 1]\\[-.05cm]
       9:&&   [&&1&& 1&& 1&& 0]\\[-.05cm]
       8:&&   [&&0&& 1&& 1&& 1]\\[-.05cm]
       7:&&   [&&1&&\;\; 0&&\;\; 0&&\;\; 1]\\[-.05cm]
       6:&&   [&&1&& 0&& 1&& 1]\\[-.05cm]
       5:&&   [&&0&& 1&& 1&& 0]\\[-.05cm]
       4:&&   [&&0&& 1&& 2&& 0]
        \end{alignat*}
        \label{tab:z222-222}
      \end{minipage}}
    \subtable[\cellquer{2}{2}{1}{3}{1}{1}]{
      \begin{minipage}[b]{.47\textwidth}
        \begin{alignat*}{6}
      11:&&\;\;    [2&& 0&& 0&& 1]\\[-.05cm]
      10:&&    [2&& 0&& 1&& 1]\\[-.05cm]
      9:&&    [1&& 1&& 1&& 0]\\[-.05cm]    
      8:&&    [1&&\;\; 1&&\;\; 2&&\;\; 0]\\[-.05cm]
      7:&&    [2&& 0&& 1&& 0]\\[-.05cm]
      6:&&    [1&& 0&& 1&& 1]\\[-.05cm]
      4:&&    [0&& 1&& 2&& 0]\\[-.05cm]
      2:&&    [1&& 0&& 2&& 0]
        \end{alignat*}
        \label{tab:z221-311}
      \end{minipage}}
  \end{minipage}
\end{table}

%%%%%%%%%%%%%%%%%%%%%%%%%%%%%%%%%%%%%%%%%%%%%%%%%%%
\subsubsection{Triangulating the Interesting Cells}
It remains to show that the pulling triangulations of
\cellquer{2}{2}{1}{3}{1}{1}, \cellquer{2}{2}{2}{2}{2}{2},
\cellquer{1}{2}{1}{2}{1}{1}, and \cellquer{2}{1}{1}{2}{1}{1} given by the
order in Table~\ref{tab:forthset} are flag. The vertex facet
incidences of these cells are listed in Table~\ref{tab:vifs} which was
generated by {\tt polymake}~\cite{polymake}.
We will repeatedly use the following fact.
\begin{lemma} \label{lem:twoFacets}
  Suppose that all but two facets $F_0$ and $F_1$ of the cell $Z$
  contain the vertex $v_0$. Pull $v_0$, and choose any triangulation
  of $Z$ refining this subdivision.
  Then every minimal non-face with more than two elements belongs to
  $F_0$ or to $F_1$. \qed
\end{lemma}

\begin{table}[hbt]
  \centering
  \caption{Vertex facet incidences of the interesting cells}
  %in Paragraph~\ref{par:problematiccells}}
  \label{tab:vifs}
  \subtable[\cellquer{1}{2}{1}{2}{1}{1}]{
    \begin{minipage}[b]{.21\textwidth}
      \centering
      \begin{tabular}{@{\extracolsep{-8pt}}
          >{$}l<{$}>{$}l<{$}>{$}l<{$}>{$}l<{$}>{$}l<{$}>{$}l<{$}}
        [v_0&v_1&v_2&v_3&v_4]\\[.05cm]
        [v_0&v_1&v_2&v_3&v_5]\\[.05cm]
        [v_0&v_1&v_4&v_5&v_6]\\[.05cm]
        [v_0&v_2&v_4&v_5&v_6]\\[.05cm]
        [v_1&v_3&v_5&\multicolumn{2}{l}{$v_6$]}\\[.05cm]
        [v_1&v_3&v_4&\multicolumn{2}{l}{$v_6$]}\\[.05cm]
        [v_2&v_3&v_4&\multicolumn{2}{l}{$v_6$]}\\[.05cm]
        [v_2&v_3&v_5&\multicolumn{2}{l}{$v_6$]}
      \end{tabular}
      \label{tab:vif-z121-211}
    \end{minipage}
    }
  \subtable[\cellquer{2}{1}{1}{2}{1}{1}]{
    \begin{minipage}[b]{.21\textwidth}
      \centering
      \begin{tabular}{@{\extracolsep{-8pt}}
          >{$}l<{$}>{$}l<{$}>{$}l<{$}>{$}l<{$}>{$}l<{$}>{$}l<{$}}
        [v_0&v_1&v_2&\multicolumn{2}{l}{$v_3$]}\\[.05cm]
        [v_0&v_1&v_3&\multicolumn{2}{l}{$v_4$]}\\[.05cm]
        [v_0&v_1&v_2&v_4&v_5]\\[.05cm]
        [v_0&v_2&v_3&\multicolumn{2}{l}{$v_6$]}\\[.05cm]
        [v_0&v_2&v_4&v_5&v_6]\\[.05cm]
        [v_0&v_3&v_4&\multicolumn{2}{l}{$v_6$]}\\[.05cm]
        [v_1&v_2&v_3&v_5&v_6]\\[.05cm]
        [v_1&v_3&v_4&v_5&v_6]
      \end{tabular}
      \label{tab:vif-z211-211}
    \end{minipage}
  }
  \subtable[\cellquer{2}{2}{2}{2}{2}{2}]{
    \begin{minipage}[b]{.23\textwidth}
      \centering
      \begin{tabular}{@{\extracolsep{-8pt}}
          >{$}l<{$}>{$}l<{$}>{$}l<{$}>{$}l<{$}>{$}l<{$}>{$}l<{$}>{$}l<{$}}
        \\[.05cm]
        [v_0&v_1&v_3&v_4&\multicolumn{2}{l}{$v_5$]}\\[.05cm]
        [v_0&v_1&v_2&v_4&\multicolumn{2}{l}{$v_6$]}\\[.05cm]
        [v_0&v_3&v_4&\multicolumn{2}{l}{$v_6$]}\\[.05cm]
        [v_0&v_1&v_2&v_3&v_5&v_7]\\[.05cm]
        [v_0&v_2&v_3&v_6&\multicolumn{2}{l}{$v_7$]}\\[.05cm]
        [v_1&v_2&v_4&v_5&v_6&v_7]\\[.05cm]
        [v_3&v_4&v_5&v_6&\multicolumn{2}{l}{$v_7$]}
      \end{tabular}
      \label{tab:vif-z222-222}
    \end{minipage}
  }
  \subtable[\cellquer{2}{2}{1}{3}{1}{1}]{
    \begin{minipage}[b]{.23\textwidth}
      \centering
      \begin{tabular}{@{\extracolsep{-8pt}}
          >{$}l<{$}>{$}l<{$}>{$}l<{$}>{$}l<{$}>{$}l<{$}>{$}l<{$}>{$}l<{$}}
        \\[.05cm]
        [v_0&v_1&v_2&v_3&\multicolumn{2}{l}{$v_4$]}\\[.05cm]
        [v_0&v_1&v_2&v_3&v_5&v_6]\\[.05cm]
        [v_0&v_1&v_4&v_5&\multicolumn{2}{l}{$v_7$]}\\[.05cm]
        [v_0&v_2&v_4&v_5&v_6&v_7]\\[.05cm]
        [v_1&v_3&v_4&\multicolumn{2}{l}{$v_7$]}\\[.05cm]
        [v_1&v_3&v_5&v_6&\multicolumn{2}{l}{$v_7$]}\\[.05cm]
        [v_2&v_3&v_4&v_6&\multicolumn{2}{l}{$v_7$]}
      \end{tabular}
      \label{tab:vif-z221-311}
    \end{minipage}
  }
\end{table}

\subsub{Triangulation of the Cell \cellquer{1}{2}{1}{2}{1}{1}} Looking
at the vertex facet incidences in Table~\ref{tab:vif-z121-211}, we see
that all  facets  opposite vertex $v_0$  are   simplices. Hence, after
pulling at  vertex $v_0$ we obtain a  simplicial complex consisting of
the $4$ simplices listed in Table~\ref{tab:triang-121211}.  Pulling at
the other vertices does not  change the complex  anymore, and its only
minimal missing faces are the edges $(v_1,v_2)$ and $(v_4,v_5)$.

\subsub{Triangulation of the Cell \cellquer{2}{1}{1}{2}{1}{1}}
The two facets opposite vertex $v_0$ are pyramids over the square
$(v_1,v_3,v_5,v_6)$ with  apex $v_2$ and $v_4$  respectively.
According to Lemma~\ref{lem:twoFacets}, any refinement of the pulling
of $v_0$ will be flag.

In our case, we obtain four facets, which are given in
Table~\ref{tab:triang-211211}. The minimal non-faces are the edges
$(v_2,v_4)$ and $(v_3,v_5)$.

\subsub{Triangulation of the Cell \cellquer{2}{2}{1}{3}{1}{1}}
The    vertex      facet  incidences     of  this     cell     are  in
Table~\ref{tab:vif-z222-222}. There are again only two facets opposite
vertex $v_0$, which are a square pyramid $S$ and a prism over a triangle
$P$. So again, after Lemma~\ref{lem:twoFacets}, we are home.

The facets of our triangulation are listed in
Table~\ref{tab:triang-221311}. The minimal non-faces are the five
edges $(v_1,v_3)$, $(v_2,v_3)$, $(v_2,v_4)$, $(v_2,v_5)$, and
$(v_5,v_6)$.

\subsub{Triangulation of the Cell   \cellquer{2}{2}{2}{2}{2}{2}}  The
vertex   facet   incidences     of   this   cell     are  in     Table
\ref{tab:vif-z221-311}.  This time  there are three facets  $S$, $P_1$
and  $P_2$ opposite  the   vertex $v_0$.  $S=(v_1,v_3,v_4,v_7)$,  is a
simplex, $P_1:=(v_1,v_3,v_5,v_6,v_7)$  is a  square pyramid with  apex
$v_7$  and $P_2:=(v_2,v_3,v_4,v_6,v_7)$ is  a square pyramid with apex
$v_3$.  Pulling at vertex  $v_0$ gives  us  a decomposition into three
cells which are pyramids over $S$, $P_1$ and  $P_2$.  The vertex $v_1$
is  contained in the base of   $P_1$, hence pulling  at $1$ decomposes
$F_1$  into two  simplices, while $F_2$  is not  affected.  Pulling at
$v_2$  decomposes $P_2$   into   two  simplices,  and  we   obtain the
simplicial complex given in Table~\ref{tab:triang-222222}.  Pulling at
the remaining  vertices  does not change   this  complex anymore. This
leaves us with the  five  minimal non-faces $(v_1,v_2)$,  $(v_2,v_5)$,
$(v_3,v_5)$, $(v_4,v_5)$ and $(v_4,v_6)$.

\begin{table}[tb]
  \caption{The  vertex  facet  incidences  of the interesting
    triangulations.}
  \label{tab:trangulations}
  \subtable[\cellquer{1}{2}{1}{2}{1}{1}]{
    \begin{minipage}[b]{.22\textwidth}
      \centering
      \begin{tabular}{@{\extracolsep{-8pt}}
         >{$}l<{$}>{$}l<{$}>{$}l<{$}>{$}l<{$}>{$}l<{$}}    
        [v_0& v_1& v_3& v_4& v_6]\\[.05cm] 
        [v_0& v_1& v_3& v_5& v_6]\\[.05cm] 
        [v_0& v_2& v_3& v_4& v_6]\\[.05cm] 
        [v_0& v_2& v_3& v_5& v_6]
      \end{tabular}
      \label{tab:triang-121211}
    \end{minipage}
  }
  \subtable[\cellquer{2}{1}{1}{2}{1}{1}]{
    \begin{minipage}[b]{.22\textwidth}
      \centering
      \begin{tabular}{@{\extracolsep{-8pt}}
          >{$}l<{$}>{$}l<{$}>{$}l<{$}>{$}l<{$}>{$}l<{$}}    
        [v_0& v_1& v_4& v_5& v_6]\\[.05cm]
        [v_0& v_1& v_3& v_4& v_6]\\[.05cm]
        [v_0& v_1& v_2& v_5& v_6]\\[.05cm]
        [v_0& v_1& v_2& v_3& v_6]
      \end{tabular}
      \label{tab:triang-211211}
    \end{minipage}
  }
  \subtable[\cellquer{2}{2}{1}{3}{1}{1}]{
    \begin{minipage}[b]{.22\textwidth}
      \centering
      \begin{tabular}{@{\extracolsep{-8pt}}
         >{$}l<{$}>{$}l<{$}>{$}l<{$}>{$}l<{$}>{$}l<{$}}    
        [v_0& v_1& v_2& v_6& v_7]\\[.05cm]
        [v_0& v_3& v_4& v_5& v_7]\\[.05cm]
        [v_0& v_3& v_4& v_6& v_7]\\[.05cm]
        [v_0& v_1& v_4& v_5& v_7]\\[.05cm]
        [v_0& v_1& v_4& v_6& v_7]
      \end{tabular}
      \label{tab:triang-221311}
    \end{minipage}
  }
  \subtable[\cellquer{2}{2}{2}{2}{2}{2}]{
    \begin{minipage}[b]{.22\textwidth}
      \centering
      \begin{tabular}{@{\extracolsep{-8pt}}
        >{$}l<{$}>{$}l<{$}>{$}l<{$}>{$}l<{$}>{$}l<{$}}          
        [v_0& v_1& v_3& v_4& v_7]\\[.05cm]
        [v_0& v_1& v_3& v_6& v_7]\\[.05cm]
        [v_0& v_1& v_5& v_6& v_7]\\[.05cm]
        [v_0& v_2& v_3& v_4& v_7]\\[.05cm]
        [v_0& v_2& v_3& v_6& v_7]
      \end{tabular}
      \label{tab:triang-222222}
    \end{minipage}
}
\end{table}

\bigskip

Hence, after pulling  at all vertices we  obtain a flag triangulation,
and, we  have  proven  the  theorem stated in  the introduction:
\begin{theorem*}[~\ref{thm:GB}]
  \textofthmquadratic\qed
\end{theorem*}
By the arguments given in the introduction, this immediately implies
\begin{theorem*}[~\ref{cor:GB}]
  \textofcorquadratic\qed
\end{theorem*}

\section{Outlook}
In some ways, these results come as  a little bit of a disappointment.
Seeing that  the toric ideal of the  Birkhoff polytope is  {\em not\/}
quadratically generated, we started this project in the hope to find a
counterexample to the conjectures  among  $3 \times 3$  transportation
polytopes. 

We know think it is conceivable to adapt the proof of
Proposition~\ref{prop:generation} to all smooth transportation
polytopes, or maybe even to general flow polytopes -- the
natural generalization of transportation polytopes.
The same seems substantially harder for the triangulation/Gr\"obner
basis result. 

In any case, the techniques can be used to improve known degree
bounds be it for sets of generators or for Gr\"obner bases: it is
sufficient to bound the degrees within the cells.

\newcommand{\etalchar}[1]{$^{#1}$}

\end{document}